\documentclass[11pt]{article}

\usepackage[dvips]{graphics}
\newcommand{\nc}{\newcommand} 
 \nc{\wh}{\widehat}
 \nc{\pl}{\partial}
 \renewcommand{\sp}{\vskip1ex}
 \nc{\inv}{^{-1}}
\def\ra{\rightarrow}
\def\iy{\infty}
\def\hf{{1\over 2}}

\def\be{\begin{equation}}
\def\ee{\end{equation}}
\def\ba{\begin{eqnarray*}}
\def\ea{\end{eqnarray*}}
\def\bae{\begin{eqnarray}}
\def\eae{\end{eqnarray}}
\def\bc{\begin{center}}
\def\ec{\end{center}}
\def\ov{\over}
\def\al{\alpha}
\def\s{\sigma}

\def\vp{\varphi}

\def\la{\lambda}
\def\pr{\textrm{Prob}}
\def\ld{\ldots}
\def\cd{\cdots}
\def\Dl{\Delta}
\def\dl{\delta}
\def\cA{\mathcal{A}}
\def\si{\sigma^{-1}}
\def\x{\xi}
\begin{document}
\title{Random Words, Toeplitz Determinants and
Integrable Systems. I.}
\date{September 27, 1999}
\author{
Alexander R.~Its\\
Department of Mathematics\\
Indiana University-Purdue University Indianapolis\\
Indianapolis, IN 46202, USA
\and
Craig A.~Tracy\\
Department of Mathematics\\
Institute of Theoretical Dynamics\\
University of California\\
Davis, CA 95616, USA
\and
Harold Widom\\
Department of Mathematics\\
University of California\\
Santa Cruz, CA 95064, USA
}
\maketitle
\begin{abstract}
It is proved that the limiting distribution of the length of the longest weakly
increasing subsequence in an inhomogeneous random word is 
related to the distribution function for the eigenvalues
of a certain \textit{direct sum} of Gaussian unitary ensembles subject
to an overall constraint that the eigenvalues lie in a 
hyperplane.
\end{abstract}
\section{Introduction}
A class of problems---important for their applications to computer science
and computational biology as well as for their inherent mathematical
interest---is the statistical analysis
of a string of random symbols.  The symbols, called \textit{letters},
are assumed to  belong to an alphabet $\cA$ of fixed size $k$.  The set of
  all such strings (or \textit{words}) of
length $N$,  $\mathcal{W}(\cA,N)$, forms the sample space in the statistical
analysis of these strings.   A natural measure on $\mathcal{W}$ is
to assign each letter equal probability, i.e.~$1/k$,  and 
to define the probability measure on words by the product measure.  Thus each letter in a word occurs
independently and with equal probability.  We call such random word models \textit{homogeneous}. 

Of course for some applications,  each letter in the alphabet does not occur with the same frequency
and  it is therefore natural to assign to each letter $i$ a probability $p_i$. If we again use the product
measure for the words (letters in a word occur  independently), then the resulting
random word models are called \textit{inhomogeneous}. 

Fixing an ordering of the alphabet $\cA$,  a \textit{weakly increasing subsequence}
of a word
 \[w=\al_1\al_2\cdots \al_N\in \mathcal{W}\]
  is a subsequence $\al_{i_1}\al_{i_2}\cdots \al_{i_m}$
such that $i_1<i_2<\cdots<i_m$ and $\al_{i_1}\le \al_{i_2}\le \cdots \le \al_{i_m}$. 
The positive integer $m$ is called the \textit{length} of this weakly increasing subsequence.
For each word $w\in\mathcal{W}$ we define $\ell_N(w)$ to equal the \textit{length of the longest
weakly increasing subsequence} in $w$.  
We  now define the fundamental object of this paper: 
\[ F_N(n):=\textrm{Prob}\left(\ell_N(w)\le n\right) \]
where $\textrm{Prob}$ is the   inhomogeneous measure on random words.  Of course,
$\textrm{Prob}$ depends  upon $N$ and the probabilities $p_i$.
 
Our results are of two types.  
To state our first results, we order the $p_i$ so that
\[ p_1\ge p_2 \ge \cdots \ge p_k \]
and decompose out alphabet $\cA$ into subsets $\cA_1$, $\cA_2$, \ldots such
that $p_i=p_j$ if and only if $i$ and $j$ belong to the same $\cA_\al$.  Setting
$k_\al=\vert\cA_\al\vert$, 
we show that the limiting distribution function as $N\ra \iy$  for
the appropriately centered and normalized random variable $\ell_N$ is related to the distribution
function for the eigenvalues $\xi_i$ in the \textit{direct sum} of mutually
independent $k_\al\times k_\al$ Gaussian unitary
ensembles (GUE),\footnote{A basic reference for random matrices is Mehta's book~\cite{mehta}.}
 conditional on the eigenvalues $\xi_i$ satisfying $\sum\sqrt{p_i}\,\xi_i=0$.
In the case when one letter occurs with greater probability than the others, this
result implies
that  the limiting distribution of $(\ell_N-N p_1)/\sqrt{N}$ is Gaussian with variance
equal to $p_1(1-p_1)$.  In the case when all the probabilities $p_i$ are distinct,
we compute the next correction in the asymptotic expansion of the mean of $\ell_N$ and find that 
\[\textrm{E}(\ell_N)=Np_1+\sum_{j>1}{p_j\ov p_1-p_j}+O({1\ov \sqrt N}),
\> N\ra\iy.\]
This last formula agrees quite well with finite $N$ simulations.  We expect this
asymptotic formula remains valid when one letter occurs with greater probability
than the others.

These results generalize work on the homogeneous model by Johansson \cite{johansson2} and
by Tracy and Widom~\cite{tw3}.  Since all the probabilities $p_i$ are equal
in the homogeneous model, the
underlying random matrix model is $k\times k$ traceless GUE.  That is,
the direct sum reduces to just one term.  In \cite{tw3} the 
integrable system underlying the finite $N$ homogeneous model was shown to be related to Painlev\'e V.  In the
isomonodromy formulation of Painlev\'e V \cite{jimbo}, the associated 
$2\times 2$ matrix linear ODE has two simple poles in the finite
complex plane and one Poincar\'e index 1  irregular singular point
at infinity.  In Part II we will show that the finite $N$ inhomogeneous model
is represented by the isomonodromy deformations of the $2\times 2$ matrix linear ODE which has 
$m+1$ simple poles in the finite complex plane and, again, one Poincar\'e index 1  
irregular singular point at infinity. The number $m$ is the total number of
the  subsets  $\cA_\al$, and the poles are located at zero point and
at the points $-p_{i_{\al}}$ ($i_{\al}$ = max $\cA_\al$). The integers
$k_{\al}$ appear as the formal monodromy exponents at the respective
points  $-p_{i_{\al}}$. We will also analyse the monodromy
meaning of the asymptotic results obtained in this part.

The results presented here are part of the recent flurry of activity centering around
connections between combinatorial probability of the Robinson-Schensted-Knuth  (RSK) type on
the one hand and random matrices and integrable systems on the other.  From the point of view
of probability theory, the quite surprising feature of these developments is that the methods
came from Toeplitz determinants, integrable differential equations of the Painlev\'e type
and the closely related Riemann-Hilbert techniques.    The first to discover
this connection at the level of distribution functions was  Baik, Deift and
Johansson \cite{bdj1} who showed that the limiting distribution of the length of the longest
increasing subsequence in a random permutation is equal to the limiting distribution function of the
appropriately centered and normalized largest eigenvalue in the GUE \cite{tw1}.  This result has
been followed by a number of 
developments relating random permutations, random words and more
generally random Young tableaux to the distribution functions
of random matrix theory \cite{bdj2, baikRains1, baikRains2, borodin2, 
grinstead, johansson1, kuperberg, okounkov, tw2}.

\setcounter{equation}{0}
\section{Random Words}
\setcounter{equation}{0}
\subsection{Probability Measure on Words and Partitions}
The Robinson-Schensted-Knuth (RSK) algorithm is a bijection between
two-line arrays $w_A$ (or
generalized permutation matrices) and ordered pairs $(P,Q)$ of
semistandard
Young tableaux (SSYT).\footnote{For a detailed account of
the RSK algorithm see Stanley, Chp.~7 \cite{stanley}.
We use without further reference various results from symmetric function theory all of
which can be found in Stanley.}  
When the two-line
arrays have the special form
\[ w_A=\left(\begin{array}{llcl}
1 & 2 & \cdots & N \\
\al_1 & \al_2 & \cdots & \al_N
\end{array}\right),\]
$\al_i\in \mathcal{A}=\{1,2,\ld,k\}$, we  identify
each $w_A$ with a word $w=\al_1\al_2\cdots\al_N$ of length $N$ composed
of letters from the alphabet $\mathcal{A}$; furthermore,
in this case  the insertion tableaux $P$ have
shape $\la\vdash N$, $\ell(\la)\le k$, with entries coming from
$\mathcal{A}$ and the
recording tableaux $Q$ are standard Young tableau (SYT) of the same
shape $\la$. As usual,  $f^\la$ denotes the number of SYT of shape $\la$
and
 $d_\la(k)$ the number of SSYT of shape $\la$ whose entries come from
$\mathcal{A}$.

We define a probability measure, $\pr$,
 on $\mathcal{W}(\mathcal{A},N)$,
 the set of all  words $w$ of length $N$
 formed from the alphabet $\mathcal{A}$, by the two requirements:
 \begin{enumerate}
 \item For each word $w$ consisting of a single letter $i\in
\mathcal{A}$,
 $\pr(w=i)=p_i$, $0<p_i<1$, with $\sum p_i=1$.
 \item For each $w=\al_1 \al_2\cdots \al_N\in\mathcal{W}$ and any
$i_j\in\mathcal{A}$, $j=1,2,\ld,N$,
\[\pr\left(\al_1 \al_2\cdots \al_N=i_1i_2\cdots i_N\right)=
\prod_{j=1}^N\pr\left(\al_j=i_j\right)\>\>\textrm{(independence)}.\]
\end{enumerate}
Of course, $\pr$ depends both on $N$ and the probabilities $\{p_i\}$.

Under the RSK correspondence, the probability measure $\pr$
induces a probability measure on partitions $\la\vdash N$, which
we will again denote by $\pr$.  This induced measure is
expressed in terms of $f^\la$ and the Schur function.  To see this we
first recall that a
tableau $T$ has \textit{type}
$\alpha=(\alpha_1,\alpha_2,\ld)$, denoted $\alpha=\textrm{type}(T)$,
if $T$ has $\alpha_i=\alpha_i(T)$ parts equal to $i$.  We write
\[ x^T=x_1^{\alpha_1(T)} x_2^{\alpha_2(T)}\cdots \]
The combinatorial definition of the Schur function of shape
$\la$ in the variables $x=(x_1,x_2,\ld)$  is the formal power series
\[ s_\la(x)=\sum_T x^T \]
summed over all SSYT of shape $\la$.  The $p=\{p_1,\ld,p_k\}$
specialization
of $s_\la(x)$ is $s_\la(p)=s_\la(p_1,p_2,\ld,p_k,0,0,\ld)$.

For each word $w\leftrightarrow (P,Q)$, the $N$ entries of $P$ consist
of the $N$ letters
of $w$ since $P$ is formed by successive row bumping the letters from
$w$.
Because of the independence assumption,
\[ p^P=p_1^{\alpha_1(P)} p_2^{\alpha_2(P)}\cdots p_k^{\alpha_k(P)} \]
gives the weight assigned to word $w$.  From the combinatorial definition
of the Schur function, we observe that its $p$ specialization is summing
the weights of  words $w$ that under RSK have shape $\lambda\vdash N$.
The recording tableau $Q$ keeps track
of the \textit{order} of the letters in the word.  The
weights of any words with the same number of letters of each type are
equal (independence), so
 we need merely count the number of such $Q$, i.e. $f^\lambda$,
 and multiply this by the weight of
any given such word to arrive at the induced measure on partitions,
\be \pr\left(\la\right)=s_\la(p) \, f^\la, \label{probPart} \ee
which satisfies the normalization $\sum_{\la\vdash N} \pr(\la)=1$.
For the homogeneous case $p_i=1/k$, the measure reduces to
\[ \pr(\la)=s_\la(1/k,1/k,\ld,1/k) \, f^\la = {d_\la(k)\, f^\la \ov
k^N}, \ \ \la\vdash N. \]
The Poissonization of this homogeneous measure is called the Charlier ensemble in
\cite{johansson2}.

If $\ell_N(w)$ equals the length of the longest \textit{weakly}
increasing subsequence in
the word $w\in \mathcal{W}(\mathcal{A},N)$, then  by the RSK
correspondence $w\leftrightarrow (P,Q)$,
the number of boxes
in the first row of $P$, $\la_1$,  equals  $\ell_N(w)$.  Hence,
\be \pr\left(\ell_N(w)\le n\right)=\sum_{{\la\vdash N \atop \la_1\le n}}
s_\la(p) \, f^\la .\label{probLength}\ee

\subsection{Toeplitz Determinant Representation}
Gessel's theorem \cite{gessel}
 is the formal power series identity\footnote{Precisely, we use the dual
version of
 Gessel's Theorem, see \S II in  \cite{tw3} whose notation we follow.}
 \[ \sum_{{\la\vdash N \atop \la_1\le n}} s_\la(x) s_\la(y) =
\det\left(T_n(\vp)\right) \]
 where $T_n(\vp)$ is the $n\times n$ Toeplitz matrix whose $i,j$ entry
is $\vp_{i-j}$, where
 $\vp_i$ is the $i^{\textrm{th}}$ Fourier coefficient of
 \[ \vp(z)=\prod_{n=1}^\iy (1+y_n z^{-1})\,\prod_{n=1}^\iy(1+x_n z), \ \
\ z=e^{i\theta}.\]

If we define the (exponential) generating function
\[ G_I(n;\{p_i\},t)=\sum_{N=0}^\infty \pr\left(\ell_N(w)\le n\right) \,
{t^N\ov N!}, \]
then an immediate consequence of Gessel's identity with $p$
specialization of the $x$ variables
and exponential specialization of the $y$ variables and the RSK
correspondence is
\be  G_I(n;\{p_i\},t)=\det\left(T_n(f_I)\right) \label{toeplitzDet}\ee
where
\be f_I(z) = e^{t/z}\, \prod_{j=1}^k (1+p_j z). \label{genFn}\ee

\section{Limiting Distribution}
\setcounter{equation}{0}
We start with the probability distribution (\ref{probPart}) on the set
of partitions
$\la=\{\la_1,\,\la_2,\ld,\la_k\}\vdash N$.  For $f^\la$ we use
the formula
\[f^\la = { N!\, \Dl(h) \ov  h_1!\,h_2! \cd h_k!}\]
where
\[ h_i=\la_j+k-i \]
and
\be \Dl(h)=\Dl(h_1,\,h_2,\ld,h_k)=\prod_{1\le i<j\le k}
(h_i-h_j).\label{Dlh}\ee
Equivalently,
\[ f^\la= {\Dl(h)\ov \prod_{i=1}^{k-1}\prod_{j=i}^{k-1}
  (\la_i+k-j)} \, {N\choose \la_1\>\>\la_2\>\>\cd \>\>\la_k} .\]

The (classical) definition of the Schur function is
\be s_{\la}(p)={\det\left(p_i^{h_j}\right)\ov \Dl(p)}
={1\ov \Dl(p)}\sum_{\s\in S_k}(-1)^\s p_1^{h_{\s(1)}}p_2^{h_{\s(2)}}\cd
p_k^{h_{\s(k)}}.
\label{schur}\ee
This holds when all the $p_i$ are distinct but in general the two
determinants require
modification, which we now describe. We order the $p_i$ so that
\be p_1\ge p_2\ge\cd \ge p_k\label{ordering}\ee
and decompose our alphabet $\cA=\{1,\,2,\ld,k\}$ into subsets
$\cA_1,\,\cA_2,\ld$
such that $p_i=p_j$ if and only if $i$ and $j$ belong to the same
$\cA_\al$. Set
$i_\al=\max\cA_\al$. Think of the $p_i$ as indeterminates and for all
indices $i$ differentiate the determinant $i_\al-i$ times with respect
to
$p_i$ if $i\in\cA_\al$. Then replace the $p_i$ by their given values.
(That
this is correct follows from l'H\^{o}pital's rule.) If we set
$k_\al=|\cA_\al|$
and write $p_\al$ for $p_{i_\al}$ then we see that $\Dl(p)$ becomes
\be\Dl'(p)=\prod_\al(1!\,2!\cd (k_\al-1)!)
\,\prod_{\al<\beta}(p_\al-p_\beta)^{k_\al\,k_\beta}\label{Dlp}\ee
and (after performing row operations) that the $i$th row of
$\det\left(p_i^{h_j}\right)$
becomes $\left(h_j^{i_\al-i}\,p_i^{h_j-i_\al+i}\right)$. Equivalently,
the partial product
$\prod_{i\in\cA_\al}p_i^{h_{\s(i)}}$
from the summand in (\ref{schur}) gets multiplied by
\be\prod_{i\in\cA_\al}\left(h_{\s(i)}^{i_\al-i}\,p_i^{-i_\al+i}\right)
=\left(\prod_{i\in\cA_\al}h_{\s(i)}^{i_\al-i}\right)
\,p_\al^{-k_\al(k_\al-1)/2}.\label{factor}\ee

In the case of distinct $p_i$ we write our formula as
\[{\rm Prob}(\la)=s_\la(p_1,\ld ,p_k)\,f^\la\]
\[={\Dl(h)\ov \Dl(p)}\,
{1\ov \prod_{i=1}^{k-1}\prod_{j=i}^{k-1}(\la_i+k-j)}\,
\sum_{\s\in S_k}(-1)^\s p_1^{k-{\s(1)}}\cdots p_k^{k-{\s(k)}}\,
p_1^{\la_{\s(1)}}\cdots p_k^{\la_{\s(k)}}\,{N\choose
\la_1\>\>\la_2\>\>\cd \>\>\la_k} .\]
Let $M_q(\la)$ denote the multinomial distribition associated with
a sequence $q=\{q_1,\,\ld,q_k\}$,
\[M_q(\la)=q_1^{\la_1}\,\cd q_k^{\la_k}\,{N\choose \la_1\>\>\la_2\>\>\cd
\>\>\la_k}.\]
If $p_\s$ denotes the sequence
$\{p_{\si(1)},\ld,p_{\si(k)}\},$ then
the above may be written
\be{\rm Prob}(\la)={\Dl(h)\ov \Dl(p)}\,{1\ov
\prod_{i=1}^{k-1}\prod_{j=i}^{k-1}(\la_i+k-j)}\,
\sum_{\s\in S_k}(-1)^\s p_1^{k-{\s(1)}}\cdots
p_k^{k-{\s(k)}}\,M_{p_\s}(\la).
\label{P}\ee
This is the formula for distinct $p_i$. In the general case we must
replace $\Dl(p)$
by $\Dl'(p)$ and each partial product
$\prod_{i\in\cA_\al}p_i^{k-{\s(i)}}$ appearing
in the sum on the right must be multiplied by the factor (\ref{factor}).

The multinomial distribution $M_q(\la)$ has the property that the total
measure
of any region where $|\la_i-N q_i|>\epsilon N$
for some $i$ and some $\epsilon>0$ tends exponentially to zero as
$N\ra\iy$. All the other terms appearing in (\ref{P}) or its
modification are uniformly bounded by a power of $N$. Since
$\la_{i+1}\le\la_i$ for all $i$ it follows that
the contribution of the terms involving $M_q(\la)$ in
(\ref{P}) will tend exponentially to zero unless $q_{i+1}\le q_i$ for
all $i$. Since
$q_i=p_{\si(i)}$ this shows that the contribution to (\ref{P}) of the
summand
corresponding to $\s$ is exponentially small unless $\s$
leaves each of the sets $\cA_\al$ invariant. It follows that if we
denote the set of such
permutations by $S'_k$ then we may restrict the sum in (\ref{P}) to the
$\s\in S'_k$ without affecting the limit. Observe that when $\s\in S'_k$
all the
$M_{p_\s}(\la)$ appearing in (\ref{P}) equal $M_p(\la)$.

Write
\[\la_i=N p_i+\sqrt{N p_i}\,\x_i.\]
In terms of the $\x_i$ the multinomial distribution $M_p(\la)$ converges
to
\be (2\pi)^{-(k-1)/2}\,e^{-\sum\x_i^2/2}\,\dl(\sum\sqrt{q_i}\,\x_i).
\label{hyperplaneGaussian}\ee
(See section~\ref{distinct}.)
Here $\dl(\sum\sqrt{q_i}\,\x_i)$ denotes Lebesgue measure on the
hyperplane $\sum\sqrt{q_i}\,\x_i=0$.

We now consider the contribution of the other terms in (\ref{P}) as
modified. Again, they are uniformly bounded by a power of $N$ and the
total measure
of any region where $|\la_i-N p_i|>\epsilon N$
for some $i$ and some $\epsilon>0$ tends exponentially to zero as
$N\ra\iy$. Thus in
determining the asymptotics of the other terms we may assume that
$\la_i\sim Np_i$
for all $i$.

The constant
$\Dl'(p)$ is given by (\ref{Dlp}). As for $\Dl(h)$, observe that
the factor
\[ h_i-h_j=\la_i-\la_j-i+j\]
in the product in (\ref{Dlh}) is asymptotically equal to $N\,(p_i-p_j)$
when $i$ and
$j$ do not belong to the same $\cA_\al$ and to
$\sqrt{N\,p_\al}\,(\x_i-\x_j)$ if
$i,\,j\in\cA_\al$. It follows that
\[\Dl(h)\sim N^{k\,(k-1)/2-\sum_\al k_\al\,(k_\al-1)/4}\,\prod
p_\al^{k_\al(k_\al-1)/4}\,
\prod_{\al<\beta}(p_\al-p_\beta)^{k_\al\,k_\beta}
\prod_\al\Dl_\al(\x),\]
where $\Dl_\al(\x)$ is the Vandermonde determinant of those $\x_i$ with
$i\in\cA_\al$.

The next factor in (\ref{P}), the reciprocal of the double product, is
asymptotically
\[N^{-k\,(k-1)/2}\,\prod_{i=1}^{k-1}p_i^{i-k}.\]

As for the sum in (\ref{P}) as modified, observe that since each $\s$
now belongs
to $S_k'$ each product appearing there is equal to $\prod p_i^{k-i}$.
Each such product is to be multiplied by
\[\prod_\al\left[\left(\prod_{i\in\cA_\al}h_{\s(i)}^{i_\al-i}\right)
\,p_\al^{-k_\al(k_\al-1)/2}\right].\]
(See (\ref{factor}).) Hence the sum itself is equal to
\[\prod_i p_i^{k-i}\,\prod_\al p_\al^{-k_\al(k_\al-1)/2}\,
\sum_{\s\in
S_k'}(-1)^\s\,\prod_\al\prod_{i\in\cA_\al}h_{\s(i)}^{i_\al-i}.\]
Since each $\s\in S_k'$ is uniquely expressible as a product of
$\s_\al\in S(\cA_\al)$
(where $S(\cA_\al)$ is the group of permutations of $\cA_\al$)
we have
\[\sum_{\s\in
S_k'}(-1)^\s\,\prod_\al\prod_{i\in\cA_\al}h_{\s(i)}^{i_\al-i}
=\prod_\al\sum_{\s_{\al}\in
S(\cA_\al)}(-1)^{\s_\al}\,\prod_{i\in\cA_\al}h_{\s_\al(i)}^{i_\al-i}\]
\[=\prod_\al \Dl_\al(h)\sim  N^{\sum k_\al(k_\al-1)/4}\,
\prod_\al\left(p_\al^{k_\al(k_\al-1)/4}\,\Dl_\al(\x)\right).\]
Putting all this together shows that the limiting distribution is
\be(2\pi)^{-(k-1)/2}\prod_\al (1!\,2!\cd (k_\al-1)!)^{-1}\,\prod_\al
\Dl_\al(\x)^2\;
e^{-\sum\x_i^2/2}\;\dl(\sum\sqrt{p_i}\,\x_i).\label{limitDistr} \ee

This has a random matrix interpretation. It is the distribution
function for the eigenvalues in the direct sum of mutually independent
$k_\al\times k_\al$
Gaussian unitary ensembles, conditional on the eigenvalues $\x_i$
satisfying
$\sum\sqrt{p_i}\,\x_i=0$.

It remains to determine the support of the limiting distribution.
In terms of the $\x_i$ the inequalities $\la_{i+1}\le\la_i$ are
equivalent to
\[\x_{i+1}\le{N(p_i-p_{i+1})\ov\sqrt Np_i}+\sqrt{p_i\ov
p_{i+1}}\,\x_{i}.\]
In the limit $N\ra\iy$ this becomes no restriction if $p_{i+1}<p_i$ but
becomes
$\x_{i+1}\le\x_i$ if $p_{i+1}=p_i$. Otherwise said, the support of the
limiting
distribution is restricted to those $\{\x_i\}$ for
which $\x_{i+1}\le\x_i$ whenever $i$ and $i+1$ belong to the same
$\cA_\al$.
(In the random matrix interpretation it means that the eigenvalues
within
each GUE are ordered.)  We denote this set of $\xi_i$ by $\Xi$.

It now follows from (\ref{probLength}) and (\ref{limitDistr}) 
(also recall the ordering (\ref{ordering})) that
\bae \lim_{N\ra\iy}\pr\left({\ell_N - N p_1 \ov \sqrt{N p_1}}\le s\right)
&=&(2\pi)^{-(k-1)/2}\prod_\al (1!\,2!\cd (k_\al-1)!)^{-1}\, \times
\label{limitingDistr} \\ 
& &\raisebox{-3mm}{$\displaystyle
\int\cdots\int \atop \hspace{-2mm} 
\textstyle{{\xi_i \in\, \raisebox{-.3ex}{$\displaystyle{\Xi}$} \atop
\vspace{.5ex} \xi_1\le s}}$} \prod_\al
\Dl_\al(\x)^2\;
e^{-\sum\x_i^2/2}\;\dl(\sum\sqrt{p_i}\,\x_i)\, d\xi_1\cdots d\xi_k 
\nonumber\eae

When the probabilities are not all equal this may be reduced to a
$k_1$-dimensional
integral as follows.
Let $i$ denote the indices in $\cA_1$ and $j$ the other indices. We have
to integrate
\[\prod_\al\Dl_\al(\x)^2\,e^{-\hf\sum \x_i^2-\hf \sum \x_j^2}\,
\dl(\sum\sqrt{p_i}\x_i+\sum\sqrt{p_j}\x_j)\]
over the subset of $\Xi$ where $\x_1\le s$. Since $\x_1=\max \x_i$ and
since the integrand 
is symmetric in 
the $\x_i$ and the $\x_j$ within their groups we may (by changing the
normalization
constant) integrate over all $\x_i\le s$ and all $\x_j$.
We first fix the $\x_i$ and integrate over the $\x_j$. These have to
satisfy
$$\sum\sqrt{p_j}\x_j=-\sum\sqrt{p_i}\x_i=-\sqrt{p_1}\sum\x_i.$$
If we write
\be\x_j=\eta_j+x\sqrt{p_j}\label{subs}\ee
where $\{\eta_j\}$ is orthogonal to $\{\sqrt{p_j}\}$ then
\[x={\sum\sqrt{p_j}\x_j\ov\sum p_j}=-{\sqrt{p_1}\ov 1-k_1 p_1}\sum
\x_i.\]
(Recall that $\cA_1$ has $k_1$ indices.) For each $\al>1$ we have
$\Dl_\al(\x)=
\Dl_\al(\eta)$ since the $p_j$ within groups are equal and
\[\sum\x_j^2=\sum\eta_j^2+x^2\sum p_j=\sum\eta_j^2+{p_1\ov 1-k_1
p_1}(\sum\x_i)^2.\]
So the distribution function is equal to a constant times
\[\int_{-\iy}^s \cd \int_{-\iy}^s\Dl(\x)^2\,
e^{-\hf[\sum \x_i^2+{p_1\ov1-k_1p_1}(\sum\x_i)^2]}\,d\x_1\cd d\x_{k_1}\,
\int\prod_{\al>1}\Dl_{\al}(\eta)^2\,e^{-\hf\sum \eta_j^2}\,d\eta,\]
where the $\eta$ integration is
over the orthogonal complement of $\{\sqrt{p_j}\}$. The $\eta$ integral
is just
another constant. Therefore the distribution function equals
\[{1\ov c_{k_1,\,p_1}}\int_{-\iy}^s \cd \int_{-\iy}^s\Dl(\x)^2\,
e^{-\hf[\sum \x_i^2+{p_1\ov1-k_1p_1}(\sum\x_i)^2]}\,d\x_1\cd
d\x_{k_1},\]
where $c_{k_1,\,p_1}$ is the integral over all of ${\bf R}^{k_1}$.

To evaluate this we make the substitution (\ref{subs}), but with $j$
replaced by $i$
and each $p_j$ replaced by $1/\sqrt k$. The integral becomes
\[\int\prod_j\Dl(\eta)^2\,e^{-\hf\sum \eta_j^2}\,d\eta\,\int
e^{-{x^2\ov2}
\left({1\ov k_1}+{p_1\ov1-k_1p_1}\right)}\,dx,\]
taken over $x\in{\bf R}$ and $\eta$ in hyperplane $\sum\eta_i=0$ with
Lebesgue measure. The $x$ integral equals $\sqrt{2\pi k_1(1-k_1p_1)}$
while the
first integral equals $(2\pi)^{(k_1-1)/2}\,1!\,2!\cd k_1!$. (For the
last, observe 
that the right side of
(\ref{limitingDistr}) must equal 1 when $s=\iy$.) Hence
\[c_{k_1,\,p_1}=(2\pi)^{k_1/2}\,1!\,2!\cd k_1!\,\sqrt{k_1(1-k_1p_1)}.\]

\subsection{\label{distinct}
Distinct probabilities---the next approximation}
If all the $p_i$ are different then $P(\la):={\rm Prob}(\la)$ equals
\be {\Dl(h)\ov \Dl(p)}\,{1\ov
\prod_{i=1}^{k-1}\prod_{j=i}^{k-1}(\la_i+k-j)}\,\prod_{i=1}^k
p_i^{k-i}\,M_p(\la)
\label{P1}\ee
plus an exponentially small correction. We recall that
\[\la_j=N p_j+\sqrt{N p_j}\,\x_j\]
and compute the Fourier transform of the measure $P$ with respect to the
$\x$ variables.
Beginning with $M_p$, we have 
\[\wh{M_p}(x)=\int e^{i\sum x_j\x_j}\,dM_p(\la)=
e^{-i\sum \sqrt{Np_j}\,x_j}\int e^{i\sum
x_j\la_j/\sqrt{Np_j}}\,dM_p(\la)\]
\[=e^{-i\sum\sqrt{Np_j}\,x_j}\,\left(\sum
p_j\,e^{ix_j/\sqrt{Np_j}}\right)^N\]
since $M_p$ is the multinomial distribution. 
An easy computation gives
\[\wh{M_p}(x)=\Big(1+{i\ov\sqrt N}Q(x)+O({1\ov N})\Big)\,
e^{-\hf\sum x_j^2+\hf(\sum\sqrt p_j x_j)^2},\]
where $Q(x)$ is a homogeneous polynomial of degree three.
(In particular the limit of $M_p$ is the inverse Fourier
transform of the exponential in the above formula,
which equals~(\ref{hyperplaneGaussian}).)

As for the other nonconstant factors in (\ref{P1}), we have 
\[\prod_{i=1}^{k-1}\prod_{j=i}^{k-1}(\la_i+k-j)=
\prod_{i=1}^{k-1}(Np_i+\sqrt{Np_i}\x_i+O(1))^{k-i}\]
\[=N^{k(k-1)/2}\prod_{i=1}^{k-1}p_i^{k-i}\,\Big(1+{1\ov\sqrt N}
\sum_{i=1}^{k-1}(k-i)\,{\x_i\ov\sqrt{p_i}}+O({1\ov N})\Big)\]
and
\[\Dl(h)=\prod_{i<j}\Big[N(p_i-p_j)+\sqrt N(\sqrt p_i\x_i-\sqrt
p_j\x_j)+O(1)\Big]\]
\[=N^{k(k-1)/2}\Dl(p)\Big(1+{1\ov\sqrt N}\sum_{i<j}{\sqrt p_i\x_i-\sqrt
p_j\x_j\ov p_i-p_j}
+O({1\ov N})\Big).\]
Thus the factors in (\ref{P1}) aside from $M_p$ contribute
\[1+{1\ov\sqrt N}\Big(\sum_{i<j}{\sqrt p_i\x_i-\sqrt p_j\x_j\ov
p_i-p_j}-
\sum_{i<j}{\x_i\ov\sqrt p_i}\Big)+O({1\ov N})
=1+{1\ov\sqrt N}\Big(\sum_{i<j}\sqrt{p_j\ov p_i}\,
{\sqrt p_j\x_i-\sqrt p_i\x_j\ov p_i-p_j}\Big)+O({1\ov N}).\]

Using the fact that multiplication by $\x_j$ corresponds, after taking
Fourier transforms, 
to $-i\pl_{x_j}$ and
combining this with the preceding we deduce that $\widehat{P}(x)$, the
Fourier 
transform of $P(\la)$ with respect to the $\x$ variables, equals 
\[\Big(1+{i\ov\sqrt N}\sum_{i<j}\sqrt{p_j\ov p_i}\,
{\sqrt p_j x_i-\sqrt p_i x_j\ov p_i-p_j}
+{i\ov\sqrt N}\,Q(x)+O({1\ov N})\Big)\,
e^{-\hf\sum x_j^2+\hf(\sum\sqrt p_j x_j)^2} \]
plus a correction which is exponentially small in $N$.\sp

\subsubsection{The mean}  
We have
\[{\rm E}(\x_1)=\int \x_1\,dP(\la)=-i\,\pl_{x_1}\widehat
P(x)\Big|_{x=0}.\]
{}From the above we see that this equals
\[{1\ov\sqrt{Np_1}}\sum_{j>1}{p_j\ov p_1-p_j}+O({1\ov N}).\]
Hence
\be \textrm{E}(\ell_N)=
{\rm E}(\la_1)=Np_1+\sum_{j>1}{p_j\ov p_1-p_j}+O({1\ov \sqrt N}),
\> N\ra\iy.\label{expectedValue}
\ee\sp
This last formula is, in fact,  an accurate approximation for 
 $\textrm{E}(\ell_N)$ (for distinct $p_i$) for moderate values of $N$.  
 Table~\ref{simulations} summarizes
 various simulations of $\ell_N$ and compares the 
 means of these simulated values with the asymptotic formula.
  We remark that even though the proof assumed distinct $p_i$,
 we expect the asymptotic formula to remain valid for $p_1>p_2\ge\cdots\ge p_k$.  
 (See the last set of simulations in Table~\ref{simulations}.)

\subsubsection{The variance}  
Let us write our approximation as
$P=P_0+N^{-1/2}P_1+O(N\inv)$ 
with corresponding expected values
${\rm E}={\rm E_0}+N^{-1/2}{\rm E_1}+O(N\inv)$. (In fact $P_1$ is a
distribution, 
not a measure, but the meaning is clear.) Then the variance of $\la_1$
is equal to
\[Np_1[{\rm E}(\x_1^2)-{\rm E}(\x_1)^2]\]
\[=Np_1\Big[{\rm E_0}(\x_1^2)-{\rm E_0}(\x_1)^2
+{1\ov\sqrt N}{\rm E_1}(\x_1^2)-{2\ov\sqrt N}{\rm E_0}(\x_1)\,{\rm
E_1}(\x_1)+O({1\ov N})
\Big].\]
Of course ${\rm E_0}(\x_1)=0$, but also
\[{\rm E_1}(\x_1^2)=-\pl^{\,2}_{x_1,\,x_1}\,\wh P_1(x)\Big|_{x=0}=0.\]
Since
\[{\rm E_0}(\x_1^2)-{\rm E_0}(\x_1)^2=1-p_1\]
we find that the variance of $\la_1$ equals $Np_1(1-p_1)+O(1)$ and so
its standard 
deviation equals $\sqrt{Np_1(1-p_1)}+O(N^{-1/2})$.

\vspace{3ex}
\textbf{\large Acknowledgments}  

This work was begun during the MSRI Semester
Random Matrix Models and Their Applications.  We wish to thank
D.~Eisenbud and H.~Rossi for their support during this semester.
This work was supported in part by the National Science Foundation
through grants DMS--9801608, DMS--9802122 and DMS--9732687.  
The last two authors thank Y.~Chen for his kind
hospitality at Imperial College where part
of this work was done as well as the  EPSRC for the award
of a Visiting Fellowship, GR/M16580, that made this visit possible.

\begin{table}
\bc
{
\begin{tabular}{|r|c|rrrr|}\hline
\raisebox{-1.5ex}[0pt]{$k$} & Probabilities  & 
\raisebox{-1.5ex}[0pt]{$N$} &  
\raisebox{-1.5ex}[0pt]{$N_S$} & 
\raisebox{-1.5ex}[0pt]{Mean}&
\raisebox{-1.5ex}[0pt]{ $\textrm{E}(\ell_N)$} \\[-1ex]
 & of \{1,\ldots,k\} & & & & \\ [0.5ex]
\hline
2 & \{ 5/7, 2/7\} & 50 & $20\,000$ & 36.37  & 36.38  \\ 
 & &		100 & $20\,000$ & 72.12 & 72.10\\
 &  &		 500 & $20\,000$ & 357.73 & 357.81 \\
 \hline
2 & \{6/11, 5/11\} & 50 & $20\,000$ & 30.54 & 32.27 \\
 & &		   100 & $20\,000$ & 58.52  & 59.55 \\
 & & 		   200 & $20\,000$ & 113.71 & 114.09 \\
 & &		400 & $20\,000$ & 223.16 & 223.18 \\
 \hline
 3 & \{1/2, 5/14, 1/7\} & 50 & $10\,000$ & 27.53 & 27.90 \\
 &                      &100 & $10\,000$ & 52.79 & 52.90 \\
 & &			500 & $10\,000$ & 252.80 & 252.90 \\
 & 	&	    1000 & $10\,000$ & 502.78    & 502.90 \\
 \hline
 3 & \{3/8, 1/3, 7/24\} & 50 & $10\,000$ & 23.96 & 30.25 \\
 & &			100&$10\,000$ & 44.33 & 49.00 \\
 & &			500& $10\,000$ & 197.65 & 199.00 \\
 & &			1000&$2\,000$ & 386.08 & 386.50 \\
 \hline
 3 & \{3/8, 5/16, 5/16\} & 50 &$10\,000$ & 23.92 & 28.75 \\
 & & 			100 & $10\,000$ & 44.16 & 47.50 \\
 & &			200 & $10\,000$ & 83.15 & 85.00 \\
 & &			400 & $10\,000$ & 159.30 & 160.00\\
 & &			800 & $10\,000$ & 310.08 & 310.00 \\
\hline
\end{tabular}
}
\ec
\caption[Simulations]{\label{simulations}
Simulations of the length of the longest weakly increasing
subsequence  in inhomogeneous random words of length $N$ for
two- and three-letter alphabets. 
$N_S$ is the sample size. The last column
gives the asymptotic  expected value (\ref{expectedValue}).}
\end{table}

\end{document}